\documentclass[aps,pre,preprint,groupedaddress]{revtex4}

\usepackage{graphicx}
\usepackage{dcolumn}
\usepackage{bm}

\begin{document}


\title{Resonant Interactions Along the Critical Line of the Riemann Zeta Function}


\author{Ronald Fisch}
\email[]{ronf124@yahoo.com}
\affiliation{382 Willowbrook Dr.\\
North Brunswick, NJ 08902}


\date{\today}

\begin{abstract}
We have studied some properties of the special Gram points of the Riemann
zeta function which lie on contour lines {\bf Im}$(\zeta ( s )) = 0$ which
do not contain zeroes of $\zeta ( s )$.   We find that certain functions of
these points, which all lie on the critical line {\bf Re}$( s )$ = 1/2, are
correlated in remarkable and unexpected ways.  We have data up to a height
of $t = 10^4$, where $s = \sigma + it$.

\end{abstract}


\maketitle

\section{Introduction}

Let $s ~=~ \sigma ~+~ it$, with $\sigma$ and $t$ real variables.
Then for $\sigma ~>~ 1$ the Riemann zeta function is defined by
\begin{equation}
  \zeta ( s ) ~=~  \sum_{n=1}^{\infty} n^{-s}  \, .
\end{equation}

It follows immediately from Eqn.~(1), that for any $t$
\begin{equation}
  \lim_{\sigma \to +\infty} \zeta ( s ) ~=~ 1  \, .
\end{equation}
It was shown by Riemann\cite{BCRW08} that $\zeta ( s )$ can be
analytically continued to a function which is meromorphic in the
complex plane, that its only divergence is a simple pole at $s$ = 1,
and that it has no zeroes on the half-plane $\sigma ~>~ 1$.

The Riemann Hypothesis\cite{BCRW08} (RH) states that the only zeroes
of $\zeta ( s )$ which do not lie on the real axis lie on the
critical line $s ~=~ {1/2} ~+~ it$.  It has resisted rigorous proof
for over 150 years, and is now widely considered to be the most
important unsolved problem in mathematics.\cite{Con03}  The
significance of the RH for physics has been shown by many
authors.\cite{BK99,ST08,MDMSWS10,SH11,FHK12}

The Gram points\cite{Edw74} are points on the critical line for which
{\bf Im}$(\zeta ( s )) = 0$ but {\bf Re}$(\zeta( s )) \ne 0$.  As stated
by H. M. Edwards,\cite{Edw74a} ``To locate the Gram points computationally
is quite easy."  The reason for this is that the spacing between neighboring
Gram points varies very smoothly, in sharp contrast to the spacing between
neighboring zeroes on the critical line.  Specifically,
\begin{equation}
  g_n - g_{n-1} ~\approx~ F ( g_{n-1} ) ~=~  {2 \pi \over {\ln ( g_{n-1} / 2 \pi )}}  \, .
\end{equation}

To show how rapidly and uniformly this converges as $n$ increases, we
display 1 minus the ratio $(g_n - g_{n-1})/F ( g_{n-1} )$ in Fig.~1.  If
we use the geometric mean $F(\sqrt{g_n g_{n-1}})$, the convergence is even
faster.  Because of this good behavior, up until now interest in the Gram
points has been focused primarily on their utility for locating the critical
zeroes of $\zeta ( s )$.

\begin{figure}
\includegraphics[width=3.4in]{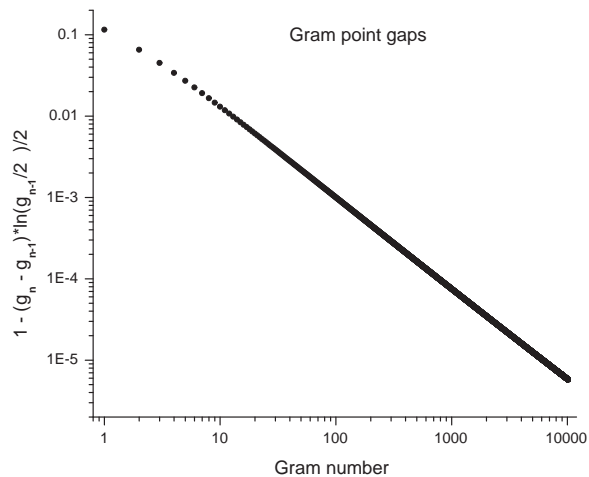}
\caption{\label{Fig.1}  1 minus the ratio between the actual gap between
two consecutive Gram points and the large $n$ approximation of Eqn.~(3),
as a function of the Gram point number, $n$.  Log-log plot.}
\end{figure}

There is a special subset of the Gram points which lie on contour lines
{\bf Im}$(\zeta ( s )) = 0$, which run from $\sigma =~- \infty$ to $\sigma =~
+ \infty$ without passing through any zero of $\zeta( s )$.\cite{Arias03,
Fisch12}  These lines divide $\zeta( s )$ into strips which run roughly
parallel to the real axis.

The contour lines {\bf Im}$(\zeta ( s )) = 0$ and {\bf Re}$(\zeta ( s )) = 0$
are highly constrained by the Cauchy-Riemann equations.  An extensive discussion
of their properties, including many illustrations, has been given by Arias de
Reyna.\cite{Arias03}  In the work presented here, the author used the computer
program of Collins to study these lines.\cite{Col09}  It is not practical to
find double precision values for the locations of the Gram points directly from
the output of Collins' computer program.  Finding approximate values from the
graphical output of this program, the accurate values used for Fig.~1 were taken
from a list supplied to the author by Michael Rubinstein.\cite{Rub13}

As shown by Arias de Reyna,\cite{Arias03} two {\bf Im}$(\zeta ( s )) = 0$
contour lines which pass through special Gram points cannot intersect when
$\sigma > 1/2$.  However, there is no general proof that contour lines of
{\bf Im}$(\zeta ( s )) = 0$ cannot intersect, although no such intersections
are known to exist.  A proof of this would immediately imply that all the
zeroes of $\zeta ( s )$ are simple, an unsettled issue which has been of
interest for a long time.  If it is possible for two {\bf Im}$(\zeta ( s )) = 0$
contour lines to intersect, then the strips bounded by the contour lines
running through the special Gram points might not be well defined.

In order for two of these {\bf Im}$(\zeta ( s )) = 0$ contour lines to
intersect at a point $s^*$, two conditions must be satisfied simultaneously.
The first condition is that $\zeta^\prime ( s^* ) = 0$, where
$\zeta^\prime ( s )$ is the derivative of $\zeta ( s )$.  The second condition
is that {\bf Im}$(\zeta ( s^* )) = 0$.  It is known\cite{LM74} that the RH
requires that there be no zeroes of $\zeta^\prime ( s )$ for $\sigma < 1/2$.
Therefore, the strips are well defined if the RH is true and in addition all
the zeroes of $\zeta ( s )$ are simple.

In this work we will study the properties of the special Gram points and their
associated contour lines.  Using data up to a height $t = 10^4$, we will find
numerical evidence of some remarkable behavior which is a result of the way the
smooth, monotonic variation in the spacings between neighboring Gram points, as
illustrated in Fig.~1, fits together with the strips bounded by the
{\bf Im}$(\zeta ( s )) = 0$ contour lines.  The average width of these strips
does not depend on the height $t$.  The author suspects that this behavior is
related in some way to Dyson's conjecture\cite{Dyson09} about the connection of
the RH with quasicrystals.

\section{Numerical results}

We define the Riemann-Siegel phase $\theta ( s )$ by
\begin{equation}
  \zeta ( s ) ~=~ |\zeta ( s )| \exp ( i \theta ( s ))  \, .
\end{equation}
Then the contour lines at the top and the bottom of each strip have
$\cos(\theta)$ = 1, {\it i.e.} $\theta$ is an integer multiple of $2 \pi$.
Counting the crossings of the critical line, $s$ = 1/2, by these contour lines
passing through the special Gram points, we plot the number of strips as a
function of the height $t$ for the first 1102 strips in Fig.~2.

\begin{figure}
\includegraphics[width=3.4in]{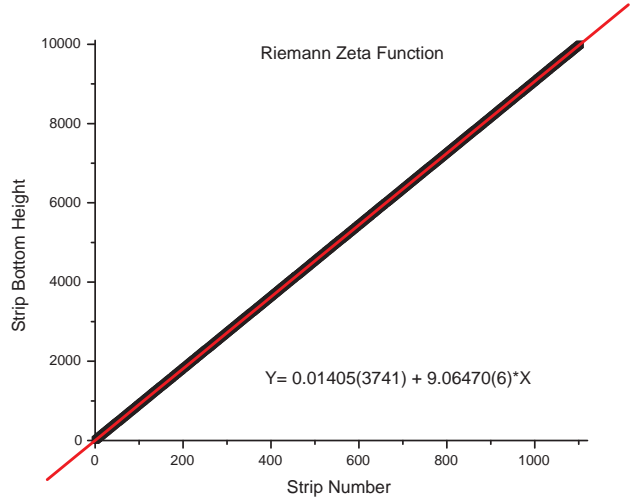}
\caption{\label{Fig.2}  Height of the bottom of a strip on the critical line
$\sigma$ = 1/2 as a function of strip number, $m$, for the first 1102 strips.}
\end{figure}

The linear least squares fit to the data is
\begin{equation}
  Y ~=~ 0.014(37) ~+~ 9.06470(6)*X  \, ,
\end{equation}
where the numbers in parentheses are the statistical errors in the last
significant figure.  Since there is a nontrivial distribution of strip
widths, there is a small amount of jitter of the data around the fitting
line.  However, there is absolutely no indication of any curvature in the
fit.  The bottom of the first strip crosses the critical line at a height
of $t$ = 9.6669080561 ... .  It is therefore somewhat mysterious that the
Y-intercept of the fitting line is consistent (within the statistical error)
with a value of zero.  If one fits the heights of the tops of the strips
instead, one finds (unsurprisingly) that the slope of the fitting line is
essentially unchanged, but the Y-intercept is now found to be 9.07(5).

For large positive $\sigma$ Eqn.~(1) can be approximated by its first two
terms.  Under this condition
\begin{equation}
  {\bf Im} (\zeta ( s )) ~\approx~  {\bf Im} (2^{-s}) ~=~ 2^{-\sigma} \sin (\ln(2)t)  \, .
\end{equation}
Thus the strip boundaries for large positive $\sigma$ and $t > 0$ will be
\begin{equation}
  t ~\approx~ 2 m \pi /\ln(2)  \, ,
\end{equation}
where the strip number, $m$, is a positive integer.  The numerical value of
$2 \pi/\ln(2)$ is 9.06472028... .  Assuming the RH is correct, it seems a
reasonable conjecture that the strips remain essentially horizontal for any
value of $t$ when $\sigma > 0$, which implies that the slope of the least
squares fit for the heights of the bottom of each strip will be independent
of $\sigma$.  The author sees no reason, however, why the Y-intercept of this
fit should be independent of $\sigma$.  In fact, it appears that for
$\sigma < 0$ this Y-intercept becomes clearly greater than zero.

The reader should note that, since the sum on the right hand side of Eqn.~(1)
does not even converge for $\sigma = 1/2$, it is rather surprising that the
data taken on the critical line, shown in Fig.~2, are well fit by a straight
line with a slope of $2 \pi / \ln (2)$.  Based on the analysis of Berry and
Keating,\cite{BK99} for example, one might have expected to see oscillations
about this line coming from the higher order terms of the sum.  Somewhat similar
ideas have been discussed by Steuding and Wegert.\cite{SW12}

We now examine in detail the departures from the straight-line behavior, to
see if there is anything resembling Berry-Keating oscillations.  We do this by
subtracting $2 m \pi/\ln(2)$ from the actual height of the $m$-th special Gram
point.  The results for various ranges of $m$ up to 1102 are shown in Fig.~3
through Fig.~7.

\begin{figure}
\includegraphics[width=3.4in]{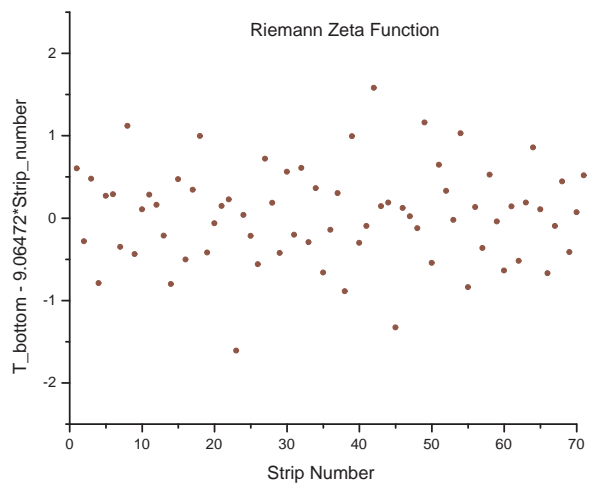}
\caption{\label{Fig.3}  Deviation of the position of the $m$-th special Gram point
from the average behavior, as a function of strip number, for the first 70 strips.}
\end{figure}

\begin{figure}
\includegraphics[width=3.4in]{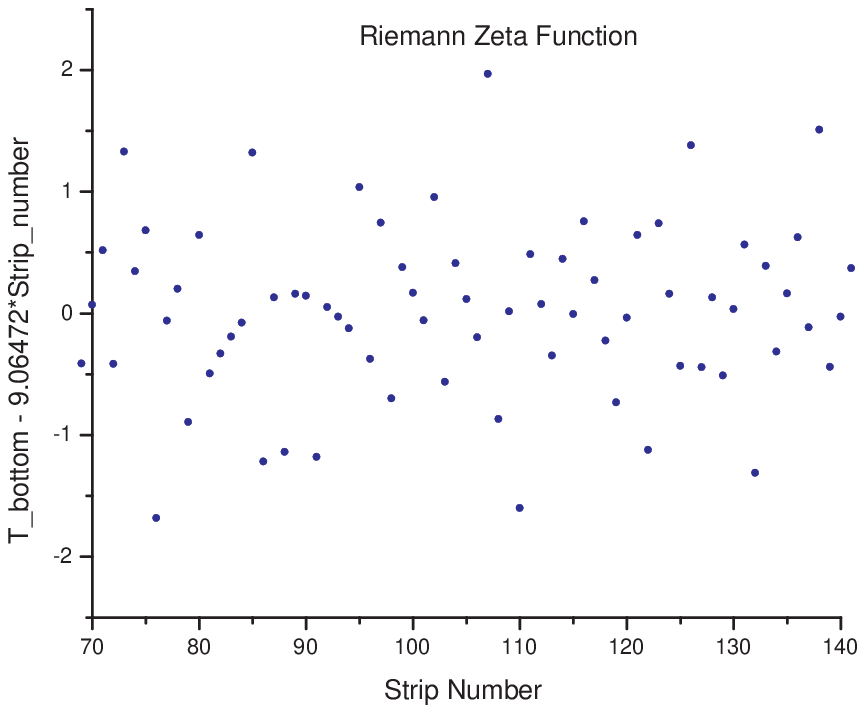}
\caption{\label{Fig.4}  Deviation of the position of the $m$-th special Gram point
from the average behavior, as a function of strip number, for strip numbers 70
to 140.}
\end{figure}

\begin{figure}
\includegraphics[width=3.4in]{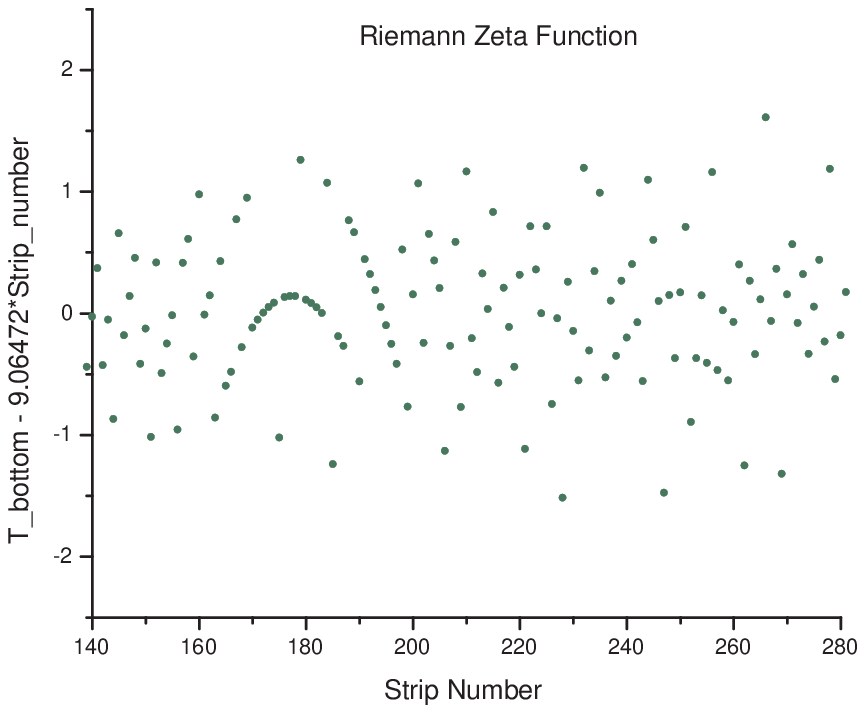}
\caption{\label{Fig.5}  Deviation of the position of the $m$-th special Gram point
from the average behavior, as a function of strip number, for strip numbers 140
to 280.}
\end{figure}

\begin{figure}
\includegraphics[width=3.4in]{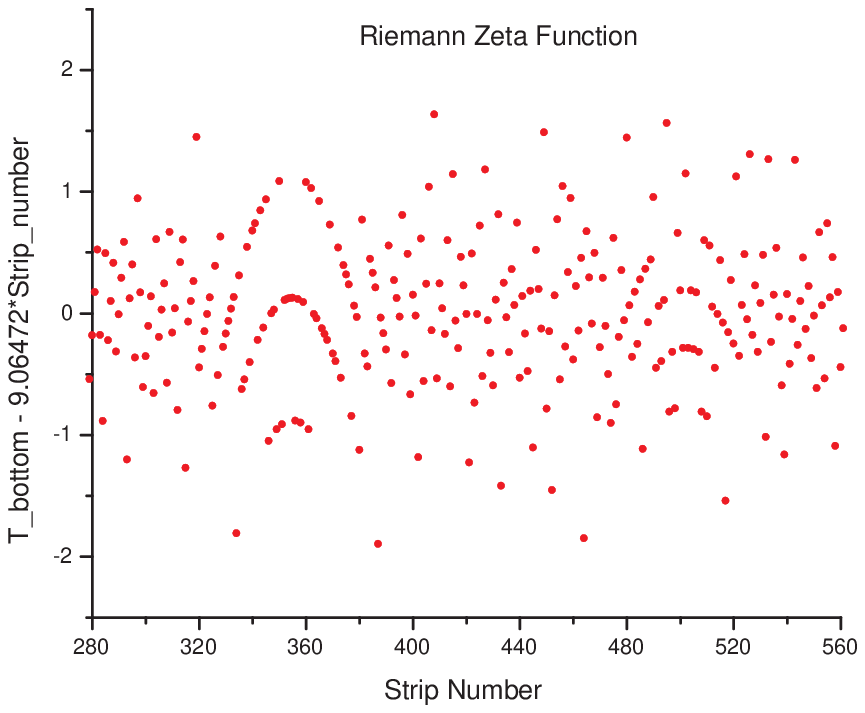}
\caption{\label{Fig.6}  Deviation of the position of the $m$-th special Gram point
from the average behavior, as a function of strip number, for strip numbers 280
to 560.}
\end{figure}

\begin{figure}
\includegraphics[width=3.4in]{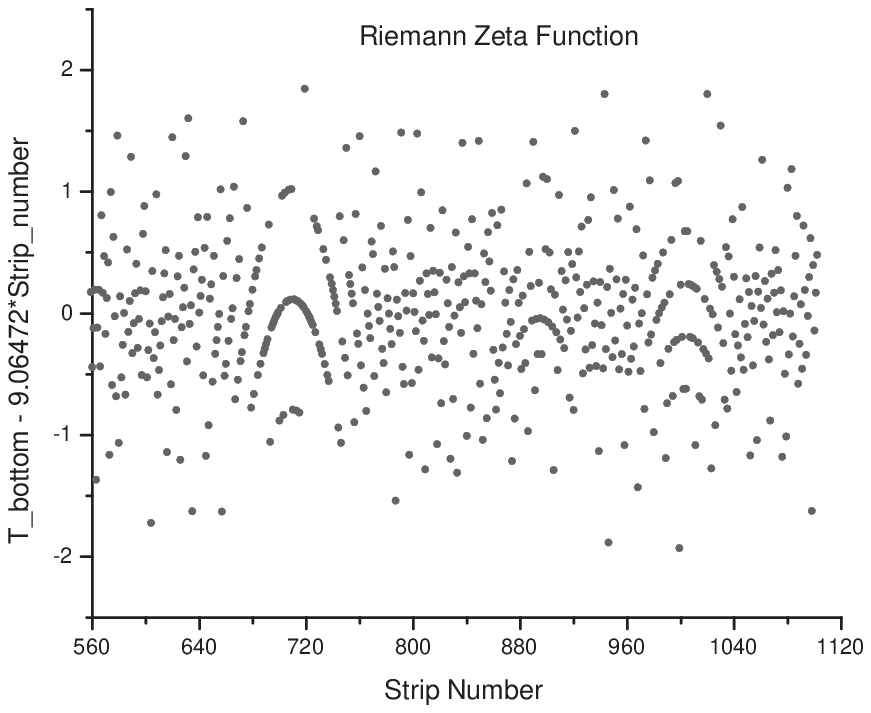}
\caption{\label{Fig.7}  Deviation of the position of the $m$-th special Gram point
from the average behavior, as a function of strip number, for strip numbers 560
to 1102.}
\end{figure}

Note that all of the points in the entire data set from 1 to 1102 have deviations
in the range -2 to 2, and that the typical size of the deviations ({\it i.e.}
the variance of the distribution of deviations) appears to be fairly insensitive
to $m$.  It is clear, however, that the points are not located randomly.  There
are obvious sets of nested arches which appear near certain values of $m$.
This demonstrates rather extensive correlations exist in those regions where
the arches are present.

It turns out that the heights where the nested arches appear are given by the
simple expression
\begin{equation}
  t_{\alpha} ~=~ 2^{1 + p/q} \pi  \, ,
\end{equation}
where $p$ and $q$ are (positive) integers with no common prime factors.  In
units of the strip number, $m$, this expression has the form
\begin{equation}
  \alpha ( p , q ) ~=~ 2^{p/q} \ln (2)  \, .
\end{equation}
The most prominent arches in Figs. 3 to 7 are centered at locations given by
$\alpha ( p , 1 )$, for $p$ = 4, 5, 6, 7, 8, 9 and 10.  Sets of secondary arches
are found at values of $\alpha ( p , 2 )$.  In Fig.~6 and Fig.~7 there are also
somewhat less defined arches which correspond to $q = 3$ and perhaps $q = 4$.
The nested arches above and below the main arches are due to strips having one
more or one less zero than they ``ought" to have.  Consistent with this idea,
the vertical spacings between the nested arches for $q = 2$ and $q = 3$ are about
1/2 and 1/3, respectively, of the vertical spacings of the $q = 1$ arches.  In
addition, all of these vertical spacings are decreasing slowly as $t$ increases.

For the $q = 1$ case Eqn.~(8) is obtained by requiring that $2 \pi /\ln(2)$
divided by the spacing between consecutive Gram points at height $t$ be equal
to $p$.  This is a kind of resonance effect between the average height of a
strip and the spacing between Gram points at height $t$.  Similarly, the
expression for larger values of $q$ may be thought of as higher order
resonances.

It would, of course, be helpful to have an explicit analytical derivation of
Eqn.~(8).  At this stage, we can only speculate about the behavior of these
arches for very large values of the height $t$.  It seems reasonable to guess
that the $q = 1$ and $q = 2$ arches will continue to be present for large
values of $t$.  It is not clear to the author what will happen at large $t$
for larger values of $q$.

\begin{figure}
\includegraphics[width=3.4in]{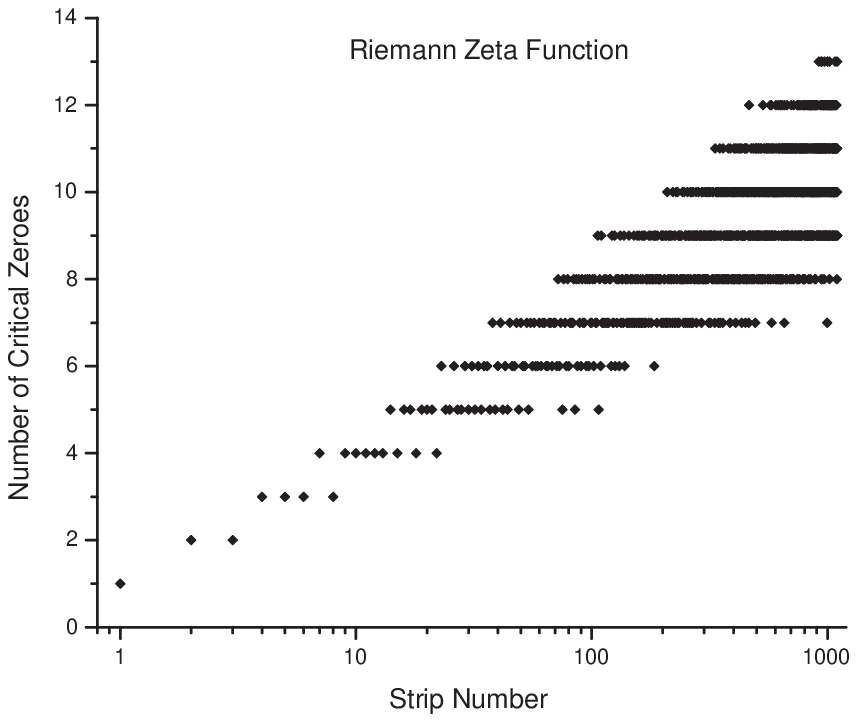}
\caption{\label{Fig.8}  Number of critical zeroes in a strip versus strip
number, for the first 1102  strips. The $X$-axis is scaled logarithmically.}
\end{figure}

The number of critical zeroes in a strip versus strip height is shown in
Fig.~8.  We see that the average number of zeroes increases logarithmically
with strip number.  It has been known for many years that the density of
the critical zeroes increases approximately as $1 / F(t)$, defined in
Eqn.~(3).  The average density of critical zeroes as a function of $t$ is
actually known to a much greater precision than this, because it is
identical to the average density of Gram points.  This follows directly
from the fact that the RH and the assumption that all zeroes are simple
implies that the number of zeroes in any strip must be equal to the number
of Gram points, counting the special Gram point on the bottom edge but
omitting the one on the top edge.

\begin{figure}
\includegraphics[width=3.4in]{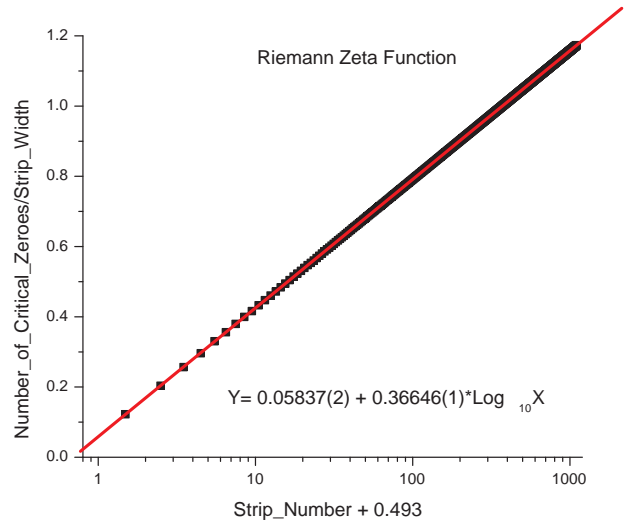}
\caption{\label{Fig.9}  (Number of zeroes on the strip)/(strip width)
versus strip number, for the first 1102 strips. The $X$-axis is scaled
logarithmically.}
\end{figure}

Due to Eqn.~(3) the width of each strip on the critical line $\sigma = 1/2$
is determined to high accuracy by the height $t$ and the number of zeroes in
the strip.  In Fig.~9 we display the function (number of zeroes on the strip)
divided by (strip width) versus the strip number.  It is thus no surprise
that these data lie on a straight line whose slope is determined by Eqn.~3,
as shown in Fig.~9.

\begin{figure}
\includegraphics[width=3.4in]{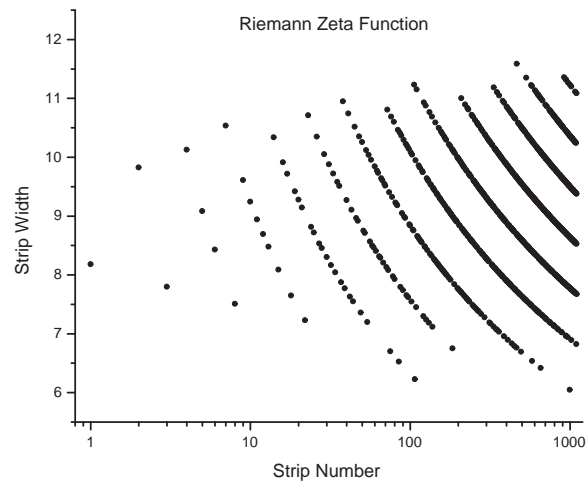}
\caption{\label{Fig.10}  Strip width on the critical line $\sigma = 1/2$ as
a function of strip number, for the first 1102 strips. The $X$-axis is scaled
logarithmically.}
\end{figure}

It is more revealing to plot the strip width versus strip number, shown in
Fig.~10.  We see that the data points sit close to lines which are hyperbolic.
Each ``line" contains all the points corresponding to a particular number of
zeroes in a strip.  The spacing between these lines is decreasing with $m$
logarithmically.  There is an observable tendency for the range of strip
widths to increase as $m$ increases, so it is not obvious what the large $m$
behavior will be.

In order to understand the behavior of the deviations of the points from the
lines, we perform an analysis similar to the one of Figs. 3 to 7.  In Fig.~11
to Fig.~15, we plot the deviations of the data from the straight line of Fig.~9.
Note that the vertical scales of Figs. 11 to 15 are all different.  The
average deviation is decreasing approximately logarithmically as $m$ increases.
This decrease compensates for the decrease in the spacing between the hyperbolic
lines as the number of zeroes per strip increases.  Therefore we expect that the
qualitative behavior of Fig.~10, ``lines" whose width is much narrower than
the spacing between them, remains valid indefinitely as $m$ increases.

\begin{figure}
\includegraphics[width=3.4in]{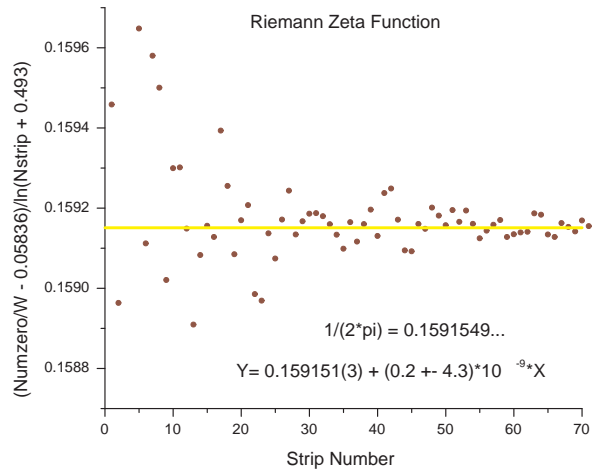}
\caption{\label{Fig.11}  Deviation of the number of zeroes on strip $m$
divided by the width of strip $m$ from a straight line fit, as a function
of strip number, for the first 70 strips.}
\end{figure}

\begin{figure}
\includegraphics[width=3.4in]{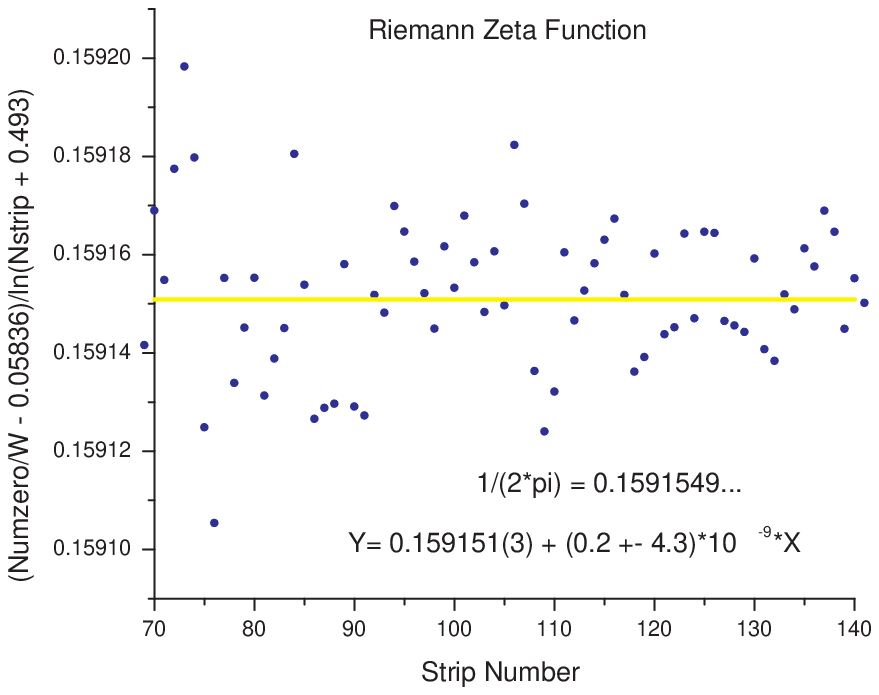}
\caption{\label{Fig.12}  Deviation of the number of zeroes on strip $m$
divided by the width of strip $m$ from a straight line fit, as a function
of strip number, for strip numbers 70 to 140.}
\end{figure}

\begin{figure}
\includegraphics[width=3.4in]{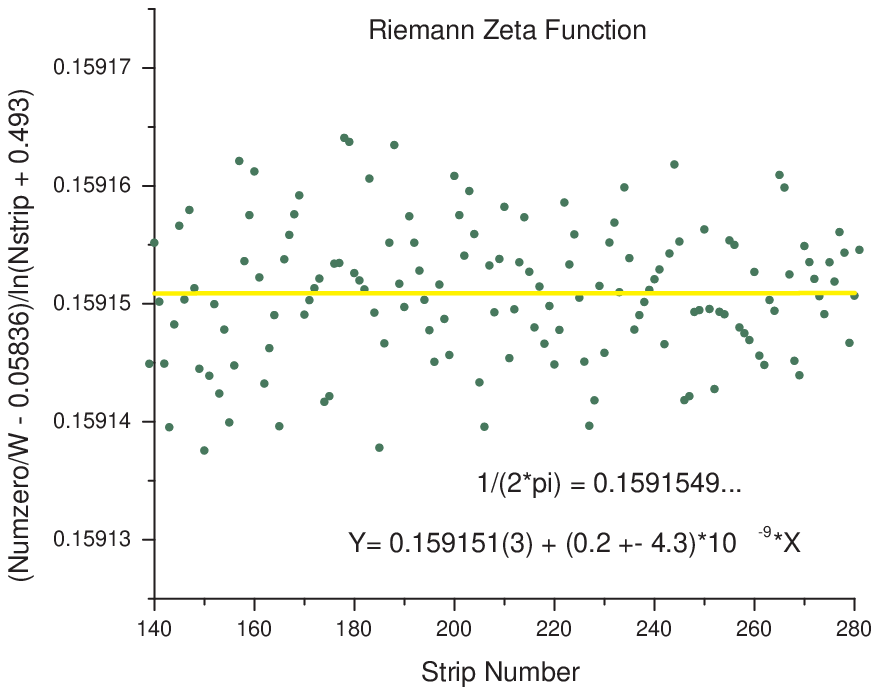}
\caption{\label{Fig.13}  Deviation of the number of zeroes on strip $m$
divided by the width of strip $m$ from a straight line fit, as a function
of strip number, for strip numbers 140 to 280.}
\end{figure}

\begin{figure}
\includegraphics[width=3.4in]{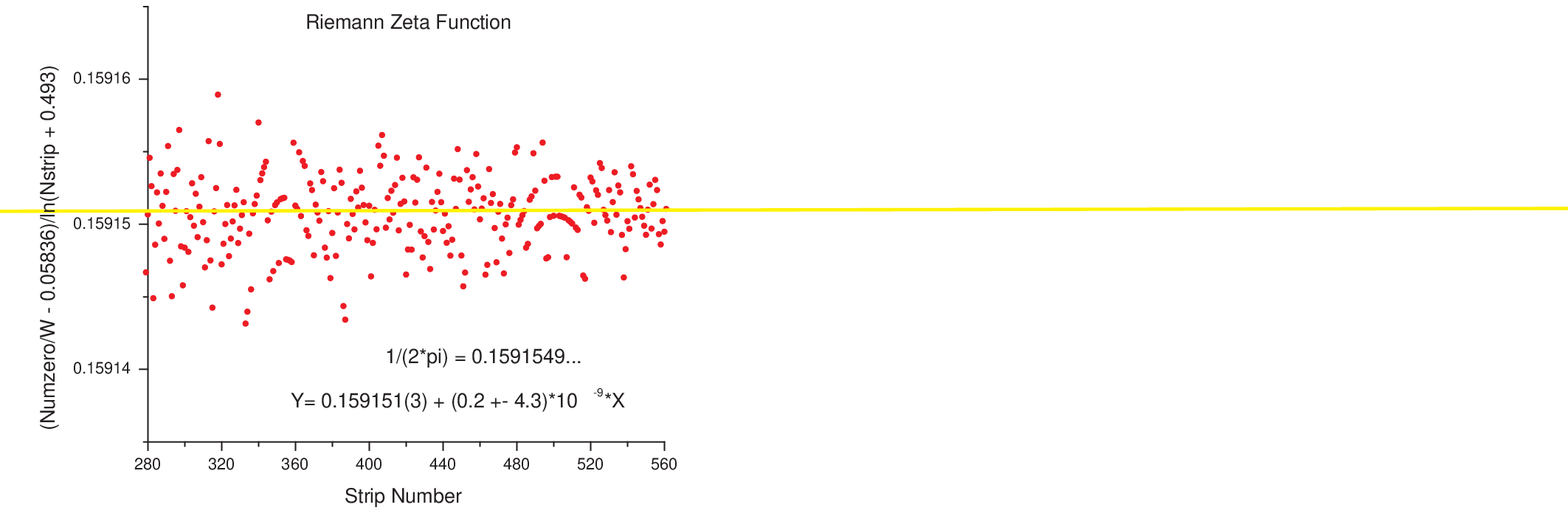}
\caption{\label{Fig.14}  Deviation of the number of zeroes on strip $m$
divided by the width of strip $m$ from a straight line fit, as a function
of strip number, for strip numbers 280 to 560.}
\end{figure}

\begin{figure}
\includegraphics[width=3.4in]{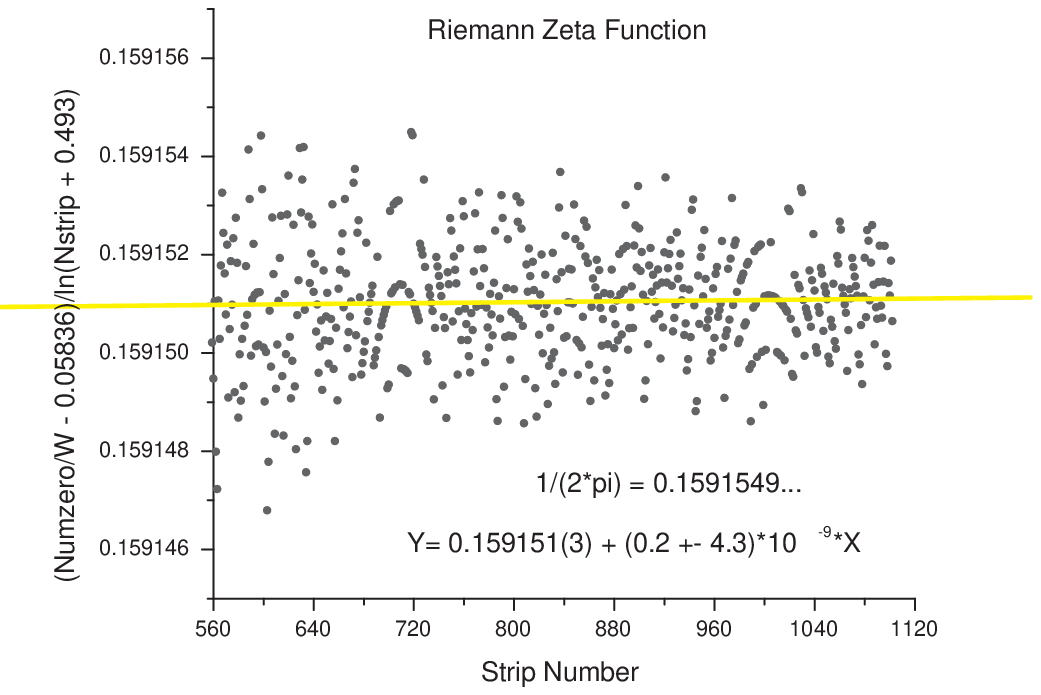}
\caption{\label{Fig.15}  Deviation of the number of zeroes on strip $m$
divided by the width of strip $m$ from a straight line fit, as a function
of strip number, for strip numbers 560 to 1102.}
\end{figure}

\begin{figure}
\includegraphics[width=3.4in]{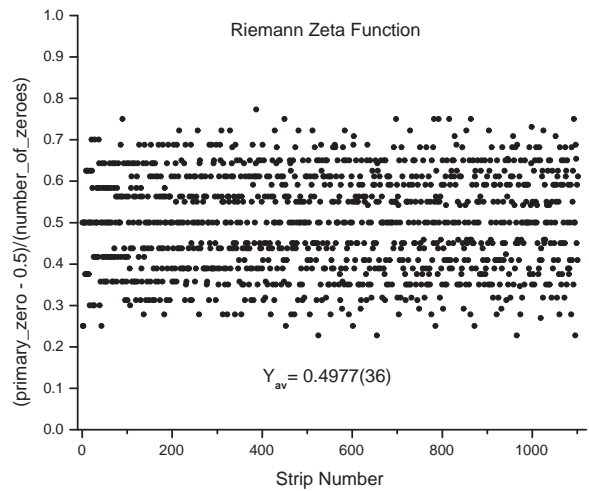}
\caption{\label{Fig.16}  (Number of the primary zero $-~0.5$)/(Number of
zeroes) versus strip number, for the first 1102 strips.}
\end{figure}

Due to the requirements of the Cauchy-Riemann equations, in each strip
there is one special zero, which we will call the primary zero.  Each
primary zero has the property that the contour with phase $\theta$ = 0
going out of it extends to $\sigma ~=~ +\infty$, at a height
\begin{equation}
  t ~\approx~ (2 m + 1) \pi /\ln(2)  \, .
\end{equation}
For all the other critical zeroes, the contour with phase $\theta ~=~ 0$
goes to $\sigma ~=~ -\infty$.

One can now ask the question ``Where is the primary zero located in the
strip?"  It seems obvious, by reason of symmetry, that the average
position of the primary zero should be at the center of the strip.
However, there is no symmetry reason why the probability distribution
for the primary zero should be uniform.  In Fig.~16 we display the values
for the function (number of the primary zero (counting from the bottom
of the strip) $-~ 0.5$) divided by (number of zeroes in the strip) versus
the strip number.  The subtraction of 0.5 in the numerator is necessary
so that this function has the value 0.5 when the primary zero is the
middle zero.

The linear least squares fit to the data of Fig.~16 shows that the average
position of the primary zero is indeed at the center of the strip.
Remarkably, one sees that the width of this probability distribution seems
to be independent of the height.  This observation is quantitatively
confirmed by calculating the variance for subsets of the data, which always
gives a result close to 0.014, independent of the range of $t$ used.

If this probability distribution remains nontrivial ({\it i.e.} neither
becoming uniform nor collapsing) in the limit $t \to \infty$, then we must
conclude that all of the zeroes in a strip are a collective entity, and
that the positions of these zeroes are highly correlated with each other.

\section{Summary}

In this work we have done an analysis of some properties of the special
Gram points of the Riemann zeta function.  We have uncovered some
previously unknown facets of the behavior of this remarkable function
along its critical line.  The most remarkable of these occur when the
ratio between $2 \pi / \ln (2)$, the average strip width, and the spacing
between Gram points passes through integers or rational numbers of small
denominator.  This seems to be some kind of a resonance effect, but its
origin is not clear at this point.

\begin{acknowledgments}
The author thanks Stephen Wolfram for making the Wolfram CDF Player 8
available as a free public download.  He also thanks Peter Sarnak and
Elliott Lieb for helpful conversations.   He thanks Michael Rubinstein
for a list of the locations of the Gram points, which made the work
described here possible.

\end{acknowledgments}



\end{document}